\documentclass{amsart}

\usepackage{amssymb,amsfonts,latexsym}

\numberwithin{equation}{section}

\newtheorem{theorem}{Theorem}[section]

\newtheorem{lemma}[theorem]{Lemma}

\newtheorem{proposition}[theorem]{Proposition}

\theoremstyle{definition}

\newtheorem{definition}[theorem]{Definition}

\newtheorem{problem}{Problem}

\newtheorem *{Theorem A}{Theorem A}
\newtheorem *{Theorem B}{Theorem B}
\newtheorem *{Corollary C}{Corollary C}
\newtheorem *{DF Conjecture}{DF Conjecture}

\newcommand{\hbt}{{\widehat {\beta}}}

\newcommand{\ben}{\begin{enumerate}}

\newcommand{\een}{\end{enumerate}}

\begin{document}

\title
{An inductive method for separable deformations}

\author{Yuval Ginosar}
\email{ginosar@math.haifa.ac.il}
\author{Ariel Amsalem}
\email{rel011235@gmail.com}
\address{Department of Mathematics, University of Haifa, Haifa 31905, Israel}
%

%


\date{\today}


\begin{abstract}
The Donald-Flanigan conjecture asserts that any group algebra of a
finite group has a separable deformation. We apply an inductive
method to deform group algebras from deformations of normal
subgroup algebras, establishing an infinite family of metacyclic
groups which fulfill the conjecture.
\end{abstract}

\maketitle
\begin{center} In memory of Murray Gerstenhaber
1927-2024.
\end{center}

\section{Introduction}

Maschke's well-known theorem asserts that the group algebra $kG$ of a finite group $G$ over a field $k$ is $k$-separable iff the order of $G$ is invertible in $k$.
In the modular case, that is where $\textnormal{char}(k)$ does divide $\left|{G}\right|$, can $kG$ still be deformed into a separable algebra?
In other words, with the notations of \S\ref{pre}, does there exist a $k[[t]]$-algebra structure on the free $k[[t]]$-module $[kG]_t:=k[[t]]\otimes_kkG$ with ${[kG]_t}/{\langle t\rangle}\cong kG$,
such that the scalar extension $\overline{k((x))}\otimes_{k[[x]]}[kG]_t$ over the algebraic closure $\overline{k((x))}$ of the field of fractions $k((x))$ of $k[[t]]$ is semi-simple?
This is known as the Donald-Flanigan [DF] problem:


\begin{DF Conjecture} (J. D. Donald and F. J. Flanigan, \cite{DF})
Let $k$ be a field and let $G$ be a finite group. Then the group algebra $kG$ admits a $k((t))$-separable deformation.
\end{DF Conjecture}
The reader is referred to the papers \cite{BG,E,ES,GG,gerstenhaber1996modular,gerstenhaber1997hecke,gerstenhaber2001donald,ginosar2019separable,peretz1999hecke,S}, which study the DF Conjecture.

This problem invites an inductive process. Given a field $k$ and an extension of groups
$1\rightarrow N \rightarrow G \rightarrow H\rightarrow 1,$
the group algebra $kG$ is isomorphic to a \textit{crossed product} $(kN)*H$.
Suppose that $kN$ satisfies the DF Conjecture. Can we make an ``inductive move"? That is, can we deduce that $(kN)*H$ also satisfies the DF Conjecture?
More specifically, let
\begin{equation}\label{exten}
1\rightarrow N \rightarrow G \rightarrow C_p\rightarrow 1
\end{equation}
be a group extension, where $C_p$ is a cyclic group of prime order $p$.
In this case, there exists an automorphism $\eta\in $Aut$_k(kN)$, where $\eta^p$ is an inner automorphism of $kN$ given as a conjugation by some $u\in (kN)^{\times}$, such that
\begin{equation}\label{int-eq1.1}
kG\cong(kN)*C_p:=kN[y;\eta]/ \left\langle y^p-u\right\rangle~
\end{equation}
is a quotient of the skew polynomial ring $kN[y;\eta]$ (see Definition \ref{pre-def2.8}). Here is a possible procedure.
\begin{problem}\label{prob}
Let $A*C_p:=A[y;\eta]/ \left\langle y^p-u\right\rangle$ be a crossed product of the cyclic group $C_p$ over a $k$-algebra $A$. Find
\begin{enumerate}
\item a separable deformation $A_t$ of $A$,
\item a deformation $\eta_t$ of $\eta$ (see Definition \ref{pre-def2.9}(2))
such that $\eta_t^p$ acts on $A_t$ as a conjugation by an element $u_t\in A_t^{\times}$ with $u_t-u\in t[A]_t$, and
\item a skew polynomial $q_t(y)\in t[A]_t[y;\eta_t]$, whose degree does not exceed $p$, turning the quotient
~$A_t[y;\eta_t]/ \left\langle y^p-u+q_t(y)\right\rangle$ into a $k((t))$-separable deformation of $A*C_p$.
\end{enumerate}
\end{problem}
This program brings some hope for a proof of the DF Conjecture for solvable groups, which are obtained by sequences of extensions of the form (\ref{exten}).
So far, this method was implemented in \cite{BG,ginosar2019separable}, where
the DF Conjecture has been proven for the family of the generalized quaternion group algebras $kQ_{2^n}\cong (kC_{2^{n-1}})*C_2$.

Once the first two steps in Problem \ref{prob} are taken care of, the third one is settled by the following theorem.
\begin{Theorem A} Let $B*C_p=B[y;\eta_t]/ \left\langle y^p-u\right\rangle$ be a crossed product, where $B$ is a $k((t))$-separable algebra.
Then there is a skew polynomial $q_t(y)\in B[y;\eta_t]$ with deg$(q_t(y))\leq 1$, such that the quotient $B[y;\eta_t]/ \left\langle y^p-u+q_t(y)\right\rangle$ is $k((t))$-separable.
Furthermore, if $B=k((t))\otimes_{k[[t]]}[B]$ for some $k[[t]]$-algebra $[B]$, then the skew polynomial $q_t(y)$ can be chosen in such a way that its coefficients belong to $t[B]$
\end{Theorem A}

Theorem A is proven in \S\ref{TA}.
Using this theorem, we are able to solve Problem \ref{prob} for $p=2$ in the case where $A:=kC_{s^2-1}$ is the group algebra of any cyclic group of the form
$C_{s^2-1}=\langle\sigma\rangle$ over every field $k$, where $u=1$ and $\eta:\sigma\mapsto\sigma^s$.
As a consequence we have

\begin{Theorem B} For any integer $s>1$, let $G_s$ be a semi-direct product $C_{s^2-1}\rtimes_\xi C_2=\left\langle\sigma\right\rangle\rtimes_\xi\left\langle\tau\right\rangle$,
where
\begin{center}
${\begin{array}{c}{\xi:{C}_{2}\mathrm{\rightarrow}{\textnormal{Aut}}\left({{C}_{s^2-1}}\right)}\\
{\mathit{\xi}\left({\mathit{\tau}}\right)\left({\mathit{\sigma}}\right)\mathrm{{=}}{\mathit{\sigma}}^{s}}\end{array}}$.
\end{center}
Then for any field $k$ of characteristic 2, the group algebra $kG_s$ has a separable deformation.
\end{Theorem B}
Theorem B is proven in \S\ref{Bproof}.
In accordance with this result, when char$(k)=2$ and $s$ is an odd integer such that $s+1$ is not a power of 2,
$\{G_s\}$ is a family of new solutions to the DF problem.

\section{Preliminaries}\label{pre}
We recall few notations and definitions, and formulate some results without proofs.
All the algebras in this note are assumed to be associative and unital. The multiplicative group of invertible elements of an algebra $A$ is denoted by $A^{\times}$.
The ring of formal power series over a field $k$ is denoted by $k[[t]]$.
Any nonzero element in the field of fractions $k((t))$ of this integral domain has the form $t^m \cdot u$ for some integer $m$ and $u\in k[[t]]^{\times}$.
Therefore, if $[A]$ is a $k[[t]]$-algebra and $x\in k((t))\otimes_{k[[t]]}[A]$, then there is an integer $m$ such that
\begin{equation}\label{m}
t^m \cdot x\in k[[t]]\otimes_{k[[t]]}[A].
\end{equation}
A (finite-dimensional) $F$-algebra $A$ is \textit{separable} over this field if the scalar extension $\overline{F}\otimes_{F}A$ of $A$ over the algebraic closure $\overline{F}$
is a semi-simple algebra.

\begin{definition}\label{pre-def2.1}
\textnormal{Let $A$ be an algebra over a field $k$. A $k((t))$-algebra $A_t$ is a \textit{deformation of A}
if there exists a $k[[t]]$-algebra structure on the free $k[[t]]$-module $[A]_t:=k[[t]]\otimes_kA$, such that
\begin{enumerate}
\item $A_t\cong k((t))\otimes_{k[[t]]}[A]_t$ as $k((t))$-algebras.
\item $A\cong [A]_t /t[A]_t$ as $k$-algebras.
\end{enumerate}
With the above notation we assume the inclusions $$A\subseteq [A]_t\subseteq A_t,$$ where $[A]_t$ is identified with $k[[t]]\otimes_{k[[t]]}[A]_t\subseteq A_t.$
In this regard, $A_t$ is called a separable deformation of $A$ if $A_t$ is a separable algebra over $k((t))$. }
\end{definition}

\begin{definition}\label{pre-def2.9}
\textnormal{With the notation of Definition \ref{pre-def2.1},\begin{enumerate}
\item A polynomial $p_t(x)\in A_t[x]$ is a deformation of $p(x)\in A[x]$ if deg$(p_t(x))=$deg$(p(x))$, and $p_t(x)-p(x)\in t[A]_t[x].$
\item
An automorphism $\eta_t\in$Aut$_{k((t))}(A_t)$ is a deformation of $\eta\in$Aut$_k(A)$ if $\eta_t(a)-\eta(a)\in t[A]_t$ for every $a\in A$.
\end{enumerate}}
\end{definition}

\begin{definition}\label{pre-def2.8}(see \cite[\S 1.2]{MR})
\textnormal{Let $R$ be a ring and let $\eta:R\to R$ be an automorphism. The set of all polynomials $f(y)$ with coefficients in $R$, endowed with the (associative) multiplication defined by
\begin{center}
$ya=\eta(a)y$, $\forall a\in R$
\end{center}
is called a \textit{skew polynomial ring} and is denoted by $R[y;\eta]$.}
\end{definition}
\begin{lemma}\label{pre-prop2.10}
With the notations of Definitions \ref{pre-def2.9} and \ref{pre-def2.8},
\begin{enumerate}
\item The algebra $A_t[x]/\left\langle p_t(x)\right\rangle$ is a deformation of $A[x]/\left\langle p(x) \right\rangle$.
\item The skew polynomial ring $A_t[y;\eta_t]$ is a deformation of $A[y;\eta]$.
\end{enumerate}
\end{lemma}

\section{Proof of Theorem A}\label{TA}
Theorem A is immediate when char$(k)\neq p$ since $B*C_p$ is already $k((t))$-separable (see \cite[Theorem 3.3]{AljadeffE.2002Ssgr}). In this case $q_t(y)=0$ does the job.
From here onwards, char$(k)= p$ is assumed.

Let $e$ be a primitive central idempotent of $B$. Then by the Wedderburn-Artin Theorem, $B\cdot e\cong M_{n_e}(D_e)$, a matrix ring
over  a division ring $D_e$. Moreover, by the separability assumption, the center $L_e:=Z(D_e)\cong Z(B\cdot e)$ is a separable field extension of $k((t))$ \cite[Theorem 6.1.2]{DK}.
Denoting the set of all primitive central idempotents of $B$ by $E$ we have
\begin{center}
$B\cong \mathop{\bigoplus}\limits_{e\in E}M_{n_e}(D_e).$
\end{center}

Clearly, $E$ is stable under the automorphism $\eta_t$.
Let $\textnormal{orb}(e)$ denote the orbit of $e\in E$ under $\eta_t$, and let
\begin{center}
$B_e:=B\cdot\sum_{e'\in \textnormal{orb}(e)}e'\cong\mathop{\bigoplus}\limits_{e'\in \textnormal{orb}(e)}M_{n_{e'}}(D_{e'})=\mathop{\bigoplus}\limits_{e'\in \textnormal{orb}(e)}M_{n_e}(D_e).$
\end{center}
Let $T\subseteq E$ be a set of representatives of the $\eta_t$-orbits, and let $u_e:=u\cdot\sum_{e'\in \textnormal{orb}(e)}e'$. Then
$B= \mathop{\bigoplus}\limits_{e\in T}B_e$, and
\begin{equation}\label{decomp}
B[y;\eta_t]/\left\langle y^p-u\right\rangle= \mathop{\bigoplus}\limits_{e\in T}B_e[y;\eta_t]/\left\langle y^p-u_e\right\rangle.
\end{equation}
The polynomial $q_t(y)$ is constructed in two steps.
In \S\ref{local} each summand in (\ref{decomp}) is handled separately, and then in \S\ref{global} a global solution is given. For each component $B_e[y;\eta_t]/\left\langle y^p-u_e\right\rangle$,
we define a polynomial $q_{t,e}(y)\in B_e[y;\eta_t]$, such that $B_e[y;\eta_t]/ \left\langle y^p +q_{t,e}(y)-u_e\right\rangle$ is $L_e$-separable.
Since $L_e$ is a separable field extension of $k((t))$, this algebra is $k((t))$-separable as well.
Afterwards, setting $q_t(y):= \sum_{e\in E}eq_{t,e}(y)$
we establish that $B[y;\eta_t]/\left\langle y^p +q_t(y)-u\right\rangle$ is separable over $k((t))$.

\subsection{The $B_e$-component.}\label{local}
Distinguish between two kinds of primitive central idempotents $e\in E$.\\
\textbf{Type 1.} Either $e$ is not stabilized by $\eta_t$ (and then $|$orb$(e)|=p$), or the restriction of $\eta_t$ to the center $Z(B_e)\cong L_e$ is non-trivial
(and then $C_p$ acts faithfully on this field).
Then by \cite[Theorem 3.3]{AljadeffE.2002Ssgr},
$B_e*C_p:=B_e[y;\eta_t]/\left\langle y^p-u_e\right\rangle$ is semi-simple over the center $Z(B_e)^{C_p}\subseteq Z(B_e)\cong L_e$.
Since this center is an intermediate field of the separable extension $k((t))\subseteq L_e$, then the summand
$B_e[y;\eta_t]/\left\langle y^p-u_e\right\rangle$ is $k((t))$-separable without any further deformation of $y^p-u_e$, that is $q_{t,e}(y)=0$.\\

\textbf{Type 2.} The automorphism $\eta_t$ both stabilizes $e$ (and so $B_e=B\cdot e\cong M_{n_e}(D_e)$) as well as acts trivially on $L_e=Z(B_e)$.
By the Skolem-Noether Theorem, $\eta_t$ is inner on the simple algebra $B_e$. Suppose that $q_e\in B_e^{\times}$ is responsible for this conjugation.
Putting $z:=q_ey$ and $v_e:=q_e^pu_e\in Z(B_e)=L_e$ we obtain
$$B_e[y;\eta_t]/ \left\langle y^p-u_e\right\rangle\cong{B_e[z]/ \left\langle z^p-v_e\right\rangle}\cong {B_e\otimes_{L_e}L_e[z]/ \left\langle z^p-v_e\right\rangle}.$$
Consider the polynomial $z^p-t^{m}z-v_e\in L_e[z]$. This polynomial is separable for any choice of positive integer $m$,
and hence, $L_e[z]/ \left\langle z^p-t^{m}z-v_e\right\rangle$ is separable over $L_e$. Tensoring it over $L_e$ with the $L_e$-separable algebra $B_e$ we obtain an $L_e$-separable,
and hence also $k((t))$-separable algebra
$$\begin{array}{rl}{B_e\otimes_{L_e}L_e[z]/ \left\langle z^p-t^{m}z-v_e\right\rangle}&\cong{B_e[z]/ \left\langle z^p-t^{m}z-v_e\right\rangle}\\
 &\cong B_e[y;\eta_t]/ \left\langle y^p-t^{m}q_e^{1-p}y-u_e\right\rangle .\end{array}$$
Thus, in this case we take $q_{t,e}(y):=-t^{m}q_e^{1-p}y$ for any positive integer $m$.

\subsection{A global solution}\label{global}

Denote the sets of idempotents in $E$ of types 1 and 2 by $E_1$ and $E_2$ respectively. By \S\ref{local}
\begin{equation}\label{qe}
q_{t,e}(y)= \left\{{\begin{array}{cc}{0}&{e\in E_1}\\{-t^{m}q_e^{1-p}y}&{e\in E_2}\end{array}}.\right.
\end{equation}

Next, let
\begin{center}
$q_t(y):=\sum_{e\in E}eq_{t,e}(y)=-\sum_{e\in E_2}et^{m}q_e^{1-p}y\in B[y;\eta_t]$
\end{center}
be a skew polynomial and note that its degree does not exceed 1. Similarly to (\ref{decomp}) we have
\begin{center}
$B[y;\eta_t]/ \left\langle y^p+q_t(y)-u\right\rangle\cong \mathop{\bigoplus}\limits_{e\in T}B_e[y;\eta_t]/ \left\langle y^p+q_{t,e}(y)-u_e\right\rangle,$
\end{center}
which by \S\ref{local} is a direct sum of $k((t))$-separable algebras and thus is $k((t))$-separable itself. The first part of the theorem is proven.

Suppose now that $B=k((t))\otimes_{k[[t]]}[B]$ for some $k[[t]]$-algebra $[B]$.
Owing to \eqref{m}, for every $e\in E$, there is a positive integer $n_e$, such that
\begin{equation}\label{ne}
t^{n_e}\cdot eq_e^{1-p}\in [B].
\end{equation}
Plug in (\ref{qe}) any positive integer $m$ satisfying
\begin{equation}\label{ta-eq10}
m> n_e,\ \ \ \forall e\in E.
\end{equation}
Then equations (\ref{qe}),(\ref{ne}) and (\ref{ta-eq10}) ensure that the coefficients of $q_t(y)$ lie in $t[B].$
 This completes the proof of Theorem A.  $\qed$

\section{Proof of Theorem B}\label{Bproof}
By (\ref{int-eq1.1}) there is an isomorphism
$kG_s\cong kC_{s^2-1}[y;\eta]/ \left\langle y^2-1\right\rangle,$ where the automorphism $\eta:kC_{s^2-1}\to kC_{s^2-1}$ is determined by the group automorphism $\xi(\tau)$, that is
$\eta(\sigma):=\sigma^s.$
Our strategy is as follows.
\begin{enumerate}
\item Deform the group $k$-algebra $kC_{s^2-1}\cong k[x]/ \langle x^{s^2-1}-1\rangle~$ to a $k((t))$-separable algebra $$(kC_{s^2-1})_t:=k((t))[x]/\langle\pi_t(x)\rangle,$$
where $\pi_t(x)\in k[[t]][x]$ is a separable polynomial of degree $s^2-1$ (see Lemma \ref{pre-prop2.10}(1)) satisfying $x^{s^2-1}-1-\pi_t(x)\in tk[[t]][x]$.
\item Construct $\eta_t\in$Aut$_{k((t))}(kC_{s^2-1})_t$ of order 2, deforming $\eta\in$Aut$_k(kC_{s^2-1})$, i.e., $\eta_t(x+\langle\pi_t(x)\rangle)=x^s+  tk[[t]][x]+ \left\langle \pi_t(x)\right\rangle$
 (see Definition \ref{pre-def2.9}(2)).
\end{enumerate}
The setup now fits that of Theorem A with $p=2$, $[B]:=[kC_{s^2-1}]_t=k[[t]][x]/\langle\pi_t(x)\rangle$ and $u:=1$.
The algebra $(kC_{s^2-1})_t[y;\eta_t]/ \left\langle y^p+1+q_t(y)\right\rangle$ established in Theorem A is a separable deformation of $kG_s$, proving Theorem B.

For establishing the two demands above we first assume that $s\neq 5$, and let
\begin{equation}\label{eq6.2}
h_t(x)=h_{t,s}(x):=x^s+tx^2.
\end{equation}
Then as can easily be verified, the polynomial
\begin{equation}\label{pit}
\pi_t(x)=\pi_{t,s}(x):=\frac{h_{t}(h_t(x))}{x}-1=\frac{(x^s+tx^2)^s+t(x^s+tx^2)^2}{x}-1
\end{equation}
has coefficients in $k[[t]]$, is of degree $s^2-1$, and defers from $x^{s^2-1}-1$ by a polynomial in $tk[[t]][x].$
Separability of $\pi_{t,s}(x)$ for $s\neq 5$ under the assumption char$(k)=2$ is established by showing that this polynomial is prime to its derivative $\pi'_t(x)$.
The proof is rather tedious and is detailed in \cite{AA}.
So much for the first step for $s\neq 5$.

Next, let $(kC_{s^2-1})_t:=k((t))[x]/\langle\pi_t(x)\rangle,$ and denote $\bar{x}:=x+\langle\pi_t(x)\rangle$. Note that by (\ref{pit})
$$h_t(\bar{x})\cdot [\bar{x}^{s-2}(\bar{x}^{s-1}+t\bar{x})^{s-1}+t(\bar{x}^{s-1}+t\bar{x})]=1,$$
in particular, $h_t(\bar{x})$ is invertible in $(kC_{s^2-1})_t$. Applying (\ref{pit}) again we have
\begin{equation}\label{htht}
h_t(h_t(\bar{x}))=\bar{x}.
\end{equation}
Define
$$\eta_t=\eta_{t,s}:
\begin{array}{ccc}
k((t))[x]/\langle\pi_t(x)\rangle&\to &k((t))[x]/\langle\pi_t(x)\rangle\\
\sum_{i=0}^{s^2-2}a_i\bar{x}^i&\mapsto &\sum_{i=0}^{s^2-2}a_ih_t(\bar{x})^i
\end{array}.$$
Owing to (\ref{htht}) and since $h_t(\bar{x})$ is invertible we have
\begin{equation}\label{van}
\pi_t(h_t(\bar{x}))=\frac{h_{t}(h_t(h_t(\bar{x})))}{h_t(\bar{x})}-1=\frac{h_{t}(\bar{x})}{h_t(\bar{x})}-1=0.
\end{equation}
The vanishing (\ref{van}) of $\pi_t(h_t(\bar{x}))$ ascertains that the well-defined map $\eta_t$ is indeed a homomorphism, which, again by (\ref{htht}), is of order 2.
Finally, since
$$\eta_t(\bar{x})-\bar{x}^s=h_t(\bar{x})-\bar{x}^s=t\bar{x}^2,$$
it follows that $\eta_t$ deforms $\eta$.
These take care of the second step, completing the proof of Theorem B for $s\neq 5$.
When $s=5$, the polynomial \eqref{pit} somehow happens to be inseparable when char$(k)=2$. Slightly modifying \eqref{eq6.2} to, e.g.,
$$h_{t,5}(x):=x^5+t(x^3+x^2),$$
and following the rest of the construction accordingly solves this problem.
 $\qed$

When does Theorem B manufacture new solutions for the DF problem?
Notice that
$$G_s\cong G_s({2'})\rtimes G_s(2),$$
where $G_s(2)$ is a semi-dihedral 2-Sylow subgroup of $G_s$, acting on the cyclic $2'$-Hall subgroup $G_s({2'})\lhd G_s$ (via its quotient of order 2) by raising its elements to their $s$-power.
\begin{itemize}
\item If $s$ is even, then the 2-Sylow subgroup of $G_s$ is cyclic. This scenario is covered in \cite{S}.
\item When $s+1$ is a power of 2, then the action of $G_s(2)$ is trivial, and hence $G_s\cong G_s({2'})\times G_s(2)$.
Tensoring the deformation for the semi-dihedral group given in \cite{E} with, say, the trivial deformation $k((t))\otimes_kkG_s({2'})$ of $kG_s({2'})$ yields a separable deformation of $kG_s$
for char$(k)$=2.
\item In the remaining cases, namely $s+1$ is even, but not a power of 2, our solutions turn out to be novel to date for char$(k)$=2.
\end{itemize}

\end{document}